\documentclass[11pt,a4paper]{article}

\usepackage[utf8]{inputenc}
\usepackage[T1]{fontenc}
\usepackage[american]{babel}

\usepackage[a4paper,margin=1in]{geometry}
\usepackage{setspace}
\usepackage{amsmath,amssymb,amsfonts,amsthm}
\usepackage{mathrsfs}
\usepackage{bbm}

\usepackage{graphicx}
\usepackage{epsfig}
\usepackage{tikz}
\usetikzlibrary{arrows,shapes,decorations.markings,decorations.pathmorphing}
\usepackage{multirow}
\usepackage{booktabs}

\usepackage[colorlinks=true,
linkcolor=blue,
citecolor=blue,
urlcolor=blue]{hyperref}

\usepackage{xcolor}
\usepackage{tikz}
\usetikzlibrary{arrows,shapes,decorations.markings,decorations.pathmorphing,positioning}

\tikzset{
	block_long/.style={
		rectangle,
		draw,
		rounded corners,
		text width=4cm,
		align=center,
		minimum height=1.2cm
	}
}

\newcommand{\bx}{\mathbf{x}}
\newcommand{\bv}{\mathbf{v}}

\newcommand{\RR}{\mathbb{R}}

\newcommand{\vf}{\varphi}

\theoremstyle{definition}

\theoremstyle{remark}

\numberwithin{equation}{section}

\begin{document}
	
	\markboth{D.~Burini, D.A.~Knopoff}{Multiscale Kinetic Structures for Living Systems}

	\title{Multiscale Kinetic Structures for Living Systems}

\author{
	Diletta Burini$^{1}$ and Damian A. Knopoff$^{2}$\\[2mm]
	$^{1}$Department of Mathematics and Computer Science, University of Perugia, Italy\\
	$^{2}$Faculty of Engineering, University of Deusto, Bilbao, Spain\\[1mm]
	\texttt{diletta.burini@unipg.it, d.knopoff@deusto.es}
}
	
	\vskip1truecm
	
	\maketitle
	

\begin{abstract}
	This paper develops a conceptual extension of the Kinetic Theory of Active Particles, building upon the framework introduced in~\cite{[BBL26]}.
	
	Living systems cannot be adequately described within classical single-scale paradigms, even when refined. To overcome this limitation, we introduce a Multiscale Kinetic Theory of Active Particles (MS-KTAP), in which a sub-microscopic scale of interacting entities is incorporated into the description of collective dynamics. In this framework, the activity variable is interpreted as an emergent quantity arising from lower-scale regulatory mechanisms and influenced by interactions across higher scales.
	
	The proposed framework captures key features of living systems—heterogeneity, adaptive decision-making, nonlinear and non-conservative interactions, spatial dynamics, and cross-scale feedback—within a unified mathematical structure. Competition and cooperation are thus described across multiple levels of organization.
	
	The first part of the paper derives the mathematical framework, while the second shows how specific models can be obtained. The paper concludes with perspectives on further developments, including possible integrations with scientific machine learning.
\end{abstract}

	\renewcommand{\arraystretch}{1.2}
	\parindent=20pt

\section{Motivations and plan of the paper} \label{Sec:1}

The Kinetic Theory of Active Particles (KTAP) is a mathematical framework developed to model the collective dynamics of living, and therefore behavioral, systems. In this setting, large populations of interacting agents are described at a mesoscopic level through distribution functions over suitable microscopic states~\cite{[BBL26]}. In this approach, interacting entities, called active particles, are characterized not only by mechanical variables such as position and velocity, but also by an internal behavioral variable, termed activity, which accounts for heterogeneous functional and social behaviors within a population.

The guiding idea of KTAP is to identify those features of complexity that are common to living systems and relevant for their collective dynamics, and to embed them into a kinetic description inspired by, yet conceptually distinct from, the classical kinetic theory of inert matter. Within this framework, KTAP has been successfully applied to a variety of contexts, including tumor--immune interactions~\cite{[BDK13]}, in-host viral dynamics~\cite{[BBO22],[BK24]}, behavioral crowds and swarms~\cite{[BGQR22]}, and behavioral economics~\cite{[BE24]}.

Despite its effectiveness, KTAP is formulated at a single mesoscopic scale, where the activity variable is treated as a primary state of the system. In biological and behavioral contexts, however, activity typically emerges as the result of complex regulatory processes acting at lower organizational levels.

Recent advances in immunology, molecular biology, and behavioral sciences indicate that functional behavior at the cellular or individual scale results from interactions occurring at finer scales.

A paradigmatic example is provided by immune checkpoints in cancer immunology. Cells of the immune system, particularly T lymphocytes, are regulated by molecular mechanisms that balance activation and inhibition signals. Proteins such as CTLA-4 and PD-1 act as negative regulators, suppressing cytotoxic responses even in the presence of antigen recognition. While these mechanisms are essential for maintaining physiological homeostasis, tumor cells can exploit them to evade immune surveillance. The seminal discoveries recognized by the 2018 Nobel Prize in Physiology or Medicine showed that targeted interventions on these regulatory pathways can restore effective immune responses~\cite{Nobel2018}.

This example highlights a fundamental limitation of single-scale descriptions: immune failure does not necessarily arise from insufficient activation at the cellular level, but may instead be the macroscopic manifestation of imbalances generated by regulatory dynamics acting at sub-cellular scales. Similar considerations apply to other classes of living systems. In the dynamics of human crowds during evacuations, for instance, individual behavior results from the interplay between internal perceptual and emotional processes and external regulatory actions, such as environmental signals or control strategies~\cite{[BGQR22]}. These processes act across different temporal and organizational scales and jointly shape the resulting collective dynamics.

These observations motivate a conceptual extension of the KTAP framework toward a genuinely multiscale perspective. Collective dynamics cannot be fully understood as the outcome of interactions occurring at a single level, but rather emerge from the coupling of processes acting across multiple interconnected scales.

Building upon the conceptual foundations developed in~\cite{[BBL26]}, the present paper introduces a Multiscale Kinetic Theory of Active Particles (MS-KTAP). In this framework, the activity variable is interpreted as an emergent quantity arising from lower-scale dynamics, while its evolution is influenced by interactions with higher-scale mechanisms. The resulting theory provides a unified mathematical setting in which competition and cooperation are not confined to interactions between populations, but extend across multiple levels of description.

The above considerations indicate that a mathematical theory of living systems cannot be obtained by a mere refinement of classical single-scale paradigms. Heterogeneity, adaptive decision-making, nonlinear and non-conservative interactions, and cross-scale regulatory feedback must therefore be incorporated at the structural level of the theory. The challenge is not simply to construct new models, but to define a mathematical architecture in which microscopic strategies, mesoscopic distributions, and macroscopic variables are consistently coupled.

A key role in this construction is played by the modeling of interactions, which represents one of the central aspects of the KTAP framework and has been critically analyzed in~\cite{[DL26]}. Within this perspective, the revisitation of Hilbert's sixth problem requires a deeper understanding of how different scales can be consistently linked within a unified mathematical structure.

After these introductory remarks, the first part of the paper is devoted to the derivation of the mathematical framework.

\vskip.1cm \noindent
Section 2 discusses the concepts of multiscale and collective dynamics in greater detail. In particular, we explain why spatial dynamics must be explicitly incorporated and why spatial homogeneity represents a restrictive assumption. The discussion is framed in the context of the interaction between mathematics and biology, drawing on contributions by Wiegner, Michael Reed~\cite{[Reed]}, Robert May~\cite{[MAY]}, and Miguel A.~Herrero~\cite{[Herrero]}.

\vskip.1cm \noindent
Section 3 presents the main result of the paper, showing how the mathematical tools of the KTAP can be extended to a multiscale setting, moving from the sub-microscopic scale to the microscopic one and leading to the derivation of kinetic-type systems. The analysis shows how the MS-KTAP structure naturally emerges from the multiscale description of interacting living systems.

\vskip.1cm \noindent
Section 4 shows how the mathematical theory developed in the previous sections can be applied to the study of self-propelled particles, such as human crowds, swarms, and programmed robotic systems. We also critically analyze the role of spatial homogeneity in specific applications. The section concludes with perspectives on future developments and applications to real-world systems.

\section{Challenges for a multiscale theory of living systems}\label{Sec:2}

This section presents the main features of the class of systems considered in this paper. The aim is to identify the conceptual framework required to derive a multiscale differential system describing the collective dynamics of living systems. These are viewed as behavioral entities and differ from classical particles in their ability to pursue objectives and to self-organize. A more in-depth philosophical discussion of these aspects can be found in pages $251$--$254$ of paragraph $65$ of the \textit{Critique of Pure Reason}, English translation~\cite{[KANT]}.

As a first step, we clarify the concept of \textit{multiscale} by referring to Hilbert's sixth problem, which poses the following two challenges:

\vskip.1cm \noindent $\bullet$ Deriving a mathematical description, possibly in the form of a differential system, that captures the dynamics across all relevant scales. This involves the microscopic (individual-based) and macroscopic (hydrodynamic) levels, with the mesoscopic scale, typical of statistical physics, acting as an intermediate bridge.

\vskip.1cm  \noindent $\bullet$ Deriving mathematical models at higher scales from the underlying dynamics at lower scales.
\vskip.1cm

The present paper mainly focuses on the first challenge. The essay by Michael Reed~\cite{[Reed]} provides useful guidelines for addressing the conceptual difficulties encountered in the development of a mathematical theory of living systems. Additional relevant contributions concern the search for causality principles~\cite{[MAY]} and the evolutionary dynamics of living systems~\cite{[MAYR]}.

These works help identify the challenges addressed in this paper, which aims to go beyond a purely biological perspective by examining the interplay between mathematics, natural sciences, and life sciences in the description of living entities. In this context, adopting a multiscale approach is essential.

Based on these considerations, we first describe the specific features that are incorporated into the mathematical framework, and then discuss the main conceptual challenges that arise in its development.

\subsection{Specific features of the system}\label{subsec:2.1}

We consider a multiscale system constituted by a large number of interacting entities. At the microscopic scale, these entities are called \textit{active particles} (a-particles), while at the lower, sub-microscopic scale, they are called \textit{sub-active particles} (sa-particles). Both types of particles may aggregate into groups that express the same function, referred to as \textit{activity}. These groups are called \textit{functional subsystems} (FSs) and \textit{sub-functional subsystems} (SFSs) at the microscopic and sub-microscopic scales, respectively. In particular, entities are grouped into FSs according to the activity they express.

The derivation of differential structures describing, within a kinetic theory framework, the collective dynamics of a-particles has been progressively developed in recent years, often motivated by applications. As discussed in Section~1, we build upon the results in~\cite{[BBL26]} to develop a multiscale differential framework including both a-particles and sa-particles, organized into FSs and SFSs.

We consider the following representations of the system:

\vskip.2cm \noindent \textit{Individual state representation:} The state of both a-particles and sa-particles is described by mechanical variables, typically position and velocity, and by behavioral variables, namely activity.

\vskip.2cm \noindent \textit{Collective state representation:} The collective state of each FS and SFS is described by a distribution function over the corresponding microscopic state space. This representation allows one to capture heterogeneity within each subsystem.
\vskip.2cm

The following features are incorporated into the mathematical description of the multiscale system:

\vskip.2cm \noindent \textit{Interactions:} Interactions occur both within and between FSs and SFSs. As a consequence, both microscopic and sub-microscopic states may be modified. These interactions may also generate competition dynamics, including proliferative and destructive events associated with mutation, selection, and post-Darwinian evolutionary mechanisms.

\vskip.2cm \noindent \textit{From interaction to decision-making:} The dynamics are driven by interactions that are modeled as learning processes, whereby each entity acquires information from neighboring entities. This learning stage is followed by a decision-making process in which actions are selected among different possible alternatives. These choices are guided by utility functions associated with the activity expressed by each entity.

\vskip.2cm \noindent \textit{The role of spatial dynamics:}
A multiscale approach requires the explicit analysis of spatial dynamics.
In particular, Giorgio Parisi's physical theories show that individual entities learn from other entities under constraints on the number of interacting neighbors. Accordingly, interactions are typically local and involve only a limited number of neighboring entities; see~\cite{[Ballerini],[Parisi2023]}.
As a consequence, learning processes are strongly influenced by spatial localization.
The simplifying assumption of spatial homogeneity must therefore be critically assessed.

\vskip.2cm \noindent \textit{Multidisciplinary aspects:}
The proposed framework applies to life sciences, social sciences, and economics. The nature of the entities, namely a-particles and sa-particles, depends on the specific system under consideration.
In biology, for instance, the microscopic scale corresponds to the cellular level, while the sub-microscopic scale corresponds to molecular dynamics.
Similarly, in the study of behavioral crowds or swarms, as well as in economic systems, the microscopic scale corresponds to individual agents, while the sub-microscopic scale corresponds to signals exchanged among individuals, such as visual, chemical, or informational cues, which influence the interaction dynamics.

A major challenge concerns multi-domain systems, where different scientific areas interact, for example physics and biology, or social dynamics and economics. In such cases, modeling must account for both cross-scale and cross-domain interactions.

\vskip.2cm \noindent \textit{Scientific machine learning:} The framework is connected with recent developments in scientific machine learning~\cite{[BLH21],[CUN24]} and builds upon ideas proposed in~\cite{[BFLT26]}. A key issue is how to incorporate learning and decision-making processes into a mathematical structure. In particular, both a-particles and sa-particles should learn the relevant variables governing the dynamics, while decision-making mechanisms should reflect their adaptive capabilities.

\subsection{Challenges and conceptual references}\label{subsec:2.2}

Beyond the specific features listed in Subsection~\ref{subsec:2.1}, the construction of a multiscale mathematical theory for living systems raises additional conceptual challenges. A first fundamental issue concerns the very notion of complexity in living systems, which lies at the core of complex systems theory. A broad perspective on this topic is offered in the work of David C.~Krakauer~\cite{[Krakauer]}, where complexity is analyzed as a structural feature emerging from the organization of interacting components across different levels.

A second key aspect is provided by Herbert A.~Simon's view of complex systems~\cite{[Simon1962],[Simon2019]}. In this perspective, complexity is not only related to the number of interacting components, but also to the hierarchical organization of their interactions and to the continuous interplay between individuals and their environment. In particular, the environment is not merely an external background, but is progressively modified by the actions and artifacts generated by the system itself. This observation is especially relevant for living systems, where feedback loops between agents and environment are essential to the emergence of collective behaviors.

These aspects are closely connected with the role of cognition. The self-organization of cognitive processes~\cite{[KA12]}, as well as their links with individual psychology~\cite{[KBOFA20]}, suggests that learning and decision-making cannot be reduced to instantaneous reaction mechanisms. Rather, they should be understood as the outcome of processes through which information is acquired, processed, and translated into adaptive responses. This viewpoint is consistent with the framework introduced in Subsection~\ref{subsec:2.1}, where interactions are interpreted not only as mechanical events, but also as learning processes followed by decision-making.

These considerations are relevant to the derivation of the differential structure underlying the proposed framework. In~\cite{[BBL26]}, learning effects are incorporated under the simplifying assumption that information is acquired rapidly, or quasi-instantaneously, so that decision-making can be directly embedded into the interaction dynamics. Although this assumption is suitable for a first kinetic description, it may become restrictive when the learning process itself evolves on a time scale that is not negligible with respect to the observed collective dynamics.

Indeed, the analysis of real systems indicates that learning may require finite time intervals, whose duration depends on the nature of the system and on the scale under consideration. This issue has been addressed in~\cite{[BFLT26]}, where the framework in~\cite{[BBL26]} is extended to account for learning processes evolving over non-negligible times. In the present setting, this point becomes even more relevant, since learning may occur simultaneously at the microscopic and sub-microscopic scales, with reciprocal influences between them.

Therefore, one of the main conceptual challenges is to embed learning, decision-making, and interaction dynamics into a genuinely multiscale structure, where different levels of organization are consistently coupled. The mathematical problem is not only to describe the dynamics at each scale, but also to model how lower-scale processes contribute to the emergence of higher-scale behaviors, and how, in turn, higher-scale dynamics feed back on the lower ones.

\section{Multiscale differential systems for collective dynamics} \label{Sec:3}

The perspective developed in Section~\ref{Sec:1} indicates that a mathematical theory of living systems must preserve heterogeneity, adaptive decision processes, nonlinear and non-conservative interactions, and feedback mechanisms across interconnected scales. These features cannot be treated as secondary modeling ingredients, but must be embedded in the structural formulation of the theory itself. We now show how these conceptual requirements lead naturally to the construction of a multiscale kinetic architecture in which microscopic strategies, mesoscopic statistical descriptions, and macroscopic observables are consistently coupled.

In this perspective, multiscale does not simply refer to the coexistence of different descriptive levels. Rather, it concerns the hierarchical organization of mechanisms acting at distinct scales, where processes at a lower level contribute to the formation of effective states at a higher level, while collective configurations emerging at higher levels influence, in turn, the evolution of lower-scale dynamics. The mathematical formulation must therefore reflect this reciprocal organization, so that cross-scale mechanisms are embedded within a coherent structural setting.

Living systems exhibit structural properties that distinguish them from classical physical systems. Individual entities are not passive units, but express strategies, adapt to environmental inputs, and modify their behavior through interactions. Heterogeneity is intrinsic: individuals may differ in their level of activation or responsiveness. Interactions are nonlinear and nonlinearly additive, so that collective outcomes cannot be reduced to the superposition of pairwise encounters. Learning mechanisms allow past interactions to influence future responses, while proliferative and destructive events continuously reshape the population. These features constitute primary sources of complexity and require mathematical structures capable of preserving heterogeneity and evolutionary dynamics.

The kinetic description provides an intermediate level between purely microscopic individual-based models and macroscopic deterministic equations. In the classical KTAP framework, interacting entities are modeled as active particles (a-particles), grouped into Functional Subsystems (FSs) indexed by $i=1,\dots,n$, each representing a population sharing a functional role. Each a-particle is characterized by a microscopic state $(\bx,\bv,u)$, where $\bx \in \Sigma \subseteq \RR^2$ denotes position, $\bv \in D_\bv$ velocity, and $u$ is an internal variable called activity. The activity encodes the strategic or functional condition of the individual and provides a mathematical representation of heterogeneity.

However, several applications reveal intrinsic limitations of a purely single-scale description. The activity variable itself may originate from regulatory mechanisms acting at lower organizational levels, such as biochemical or informational processes. Moreover, interaction rules at the microscopic level may depend on signals produced by the population or on macroscopic observables. These considerations motivate a genuinely multiscale extension.

We therefore consider a system operating at two structurally distinct yet dynamically coupled scales.

At the microscopic scale (m-scale), the interacting entities are the a-particles described above. Their state is given by $(\bx,\bv,u)$, and they are grouped into FSs indexed by $i=1,\dots,n$.

At a lower level, called sub-microscopic scale (sm-scale), we introduce sub-active particles (sa-particles), grouped into Sub-Functional Subsystems (SFSs) indexed by $j=1,\dots,m$. These entities represent regulatory signals produced by a-particles and capable of influencing interaction rules at the microscopic level. Each sa-particle is characterized by a state $(\bx,\bv,w)$, where $w$ measures the intensity or effectiveness of the signal.

The activity variables $u$ and $w$ are assumed to take values in bounded intervals. If $u_m$ and $u_M$ denote the minimum and maximum admissible values of the microscopic activity, and $w_m$ and $w_M$ the corresponding bounds for the sub-microscopic activity, we introduce the normalized variables
\begin{equation}\label{normalized_variables}
	\hat{u} = \frac{u - u_m}{u_M - u_m},
	\qquad
	\hat{w} = \frac{w - w_m}{w_M - w_m},
\end{equation}
so that $\hat{u},\hat{w} \in [0,1]$. For notational simplicity, we continue denoting the normalized variables by $u$ and $w$. This normalization allows both activity variables to be treated as dimensionless quantities while preserving their interpretation as structured internal states defined at different organizational levels.

The m-scale and sm-scale represent different levels of organization within the same living system. They are dynamically coupled: a-particles generate sa-particles, while sa-particles influence the evolution of a-particles through modified interaction rules.

This coupling cannot be interpreted as a simple superposition of effects. The evolution of the sm-scale is triggered by interaction processes occurring at the microscopic level, while the distribution of sa-particles progressively modifies the effective interaction structure experienced by the a-particles. In this way, the two scales do not operate independently: the regulatory signals generated at the lower level become part of the interaction environment at the higher level, thereby contributing to shape subsequent decision processes. The resulting organization reflects a closed dynamical structure in which production and regulation are intrinsically intertwined.

The statistical state of the system is described by distribution functions over the corresponding phase spaces. For the microscopic scale, we introduce
\begin{equation}\label{f_i}
	f_i = f_i(t,\bx,\bv,u),
	\qquad i=1,\dots,n,
\end{equation}
while for the sub-microscopic scale we define
\begin{equation}\label{varphi_j}
	\varphi_j = \varphi_j(t,\bx,\bv,w),
	\qquad j=1,\dots,m.
\end{equation}
These functions represent the density of entities in the extended microscopic state spaces.

Macroscopic quantities are obtained as moments of the distribution functions. For instance, the local density of the $i$-th FS is given by
\begin{equation}\label{rho_i}
	\rho_i(t,\bx)
	=
	\int_{D_\bv \times D_u}
	f_i(t,\bx,\bv,u)\, d\bv\,du,
\end{equation}
and analogous expressions hold for the SFSs. Mean velocity and mean activity fields are derived by suitable weighted moments.

Within this framework, interactions are interpreted as stochastic decision processes inspired by game-theoretical concepts. When a-particles encounter other individuals or perceive environmental signals within their interaction domain, they update their microscopic state according to probabilistic rules. The outcome of interactions is modeled through interaction rates and transition probabilities, which may depend on the microscopic states of the interacting entities as well as on the evolving distribution functions. Such interactions may induce changes in activity, transitions across subsystems, and proliferative or destructive events. Learning and selection mechanisms are thus embedded directly into the kinetic description.

A key conceptual element inherited from KTAP is the separation between learning and decision-making: entities first acquire information within their interaction domain and subsequently update their state according to decision rules that may evolve in time. In the multiscale setting, this mechanism is influenced by the presence of sub-microscopic sa-particles, which modify the effective interaction operators.

In order to formalize interactions within the KTAP framework, it is necessary to distinguish among different roles played by interacting entities. At the microscopic scale, three types of a-particles are considered.

\textit{Test a-particles}, described by the distribution function $f_i(t,\bx,\bv,u)$, are representative of the whole system, with the index $i$ identifying the corresponding functional subsystem.

\textit{Candidate a-particles}, described by $f_i(t,\bx_\ast,\bv_\ast,u_\ast)$, are the particles that may acquire, in probability, the micro-state of the test particles after interaction with the field particles.

Finally, \textit{field a-particles}, described by $f_i(t,\bx^\ast,\bv^\ast,u^\ast)$, interact with the test particles and may lose their micro-state as a consequence of the encounter. The identification of candidate and field particles is performed in probability, through their distribution over the microscopic state.

A similar structure is introduced at the sub-microscopic scale. Although interactions of sa-particles may differ in nature from those of a-particles, for the purposes of a general multiscale theory they can be described according to the same interaction paradigm.

Accordingly, \textit{test sa-particles} are represented by $\varphi_j(t,\bx,\bv,w)$, with $j$ denoting the sub-functional subsystem. \textit{Candidate sa-particles}, described by $\varphi_j(t,\bx_\ast,\bv_\ast,w_\ast)$, may acquire, in probability, the micro-state of the test sa-particles after interaction with field sa-particles.

\textit{Field sa-particles}, described by $\varphi_j(t,\bx^\ast,\bv^\ast,w^\ast)$, interact with the test sa-particles and may lose their micro-state, their identification being determined probabilistically by their distribution over the microscopic state.

The interaction mechanisms described above affect not only the microscopic state of the entities, but also the structure of the statistical flows in the extended phase space. Microscopic encounters may generate regulatory signals at the sub-microscopic level, while the evolving distribution of such signals modifies the effective transition probabilities governing subsequent interactions. As a consequence, the balance of entities at each scale must account simultaneously for intra-scale interactions and for cross-scale exchange mechanisms, including processes that are not necessarily conservative due to proliferation, destruction, or regulatory modulation.

\vskip.2cm

In order to specify the structure of interactions, we assume that each a-particle and each sa-particle is endowed with a \textbf{sensitivity domain}.

At the microscopic scale (FS), the sensitivity domain is denoted by
\begin{equation}\label{omega}
	\Omega = \Omega(\bx,\omega),
\end{equation}
where $\omega$ denotes the velocity direction of the a-particle.

At the sub-microscopic scale (SFS), the sensitivity domain is denoted by
\begin{equation}\label{omega_m}
	\Omega_m = \Omega_m(\bx,\omega_m),
\end{equation}
where $\omega_m$ denotes the velocity direction of the sa-particle.

These domains depend on the position $\bx$ and on the corresponding velocity direction, and represent the region in which the entity is able to perceive and interact with other entities.

In spatially distributed systems, $\Omega$ may be modeled as a sector (in two space dimensions) or a cone (in three space dimensions) with vertex at $\bx$ and axis aligned with the direction $\omega$. It therefore defines the region in which information is acquired before a state update occurs.

Interactions are therefore nonlocal: the evolution of the state at $(t,\bx,\bv,u)$ may depend on the distribution of particles and sa-particles in a neighborhood determined by $\Omega$. Moreover, unlike classical mechanical collisions, interactions are in general not reversible and may depend not only on the states of the interacting pair, but also on the local statistical configuration of the system, that is, on the distribution functions themselves.

Within this framework, three main classes of interactions are considered:

\noindent
\textit{Microscopic interactions within and across functional subsystems:}
interactions may occur among a-particles belonging to the same FS or to different FSs.

\noindent
\textit{Sub-microscopic interactions within and across sub-functional subsystems:}
sa-particles may interact among themselves, either within the same SFS or across different SFSs.

\noindent
\textit{Cross-scale interactions:}
a-particles and sa-particles may interact across scales. In particular, a-particles may generate sa-particles, while the distribution of sa-particles may influence the interaction rules governing a-particles.

From the viewpoint of their structural effects on the system, interactions can be further classified into two fundamental types.

\vskip.2cm \noindent \textit{Number-conservative interactions.}
These interactions modify the microscopic state --- namely velocity and activity --- of the interacting entities without altering their total number. They account for behavioral changes, strategy adaptation, and learning processes, while preserving the population size within each functional subsystem.

\vskip.2cm \noindent \textit{Proliferative and destructive interactions.}
These interactions lead to a gain or loss of entities, possibly following encounters across different subsystems. They model birth, death, mutation, or suppression phenomena, which are intrinsic to living systems and are generally not mass-conservative. In the multiscale setting, such mechanisms may act both at the microscopic level (a-particles) and at the sub-microscopic level (sa-particles), depending on the specific application under consideration.

\vskip.2cm

A key structural feature of the interaction dynamics is the separation between \emph{learning} and \emph{decision-making}. When an entity enters its sensitivity domain $\Omega$, it first acquires information about the microscopic and possibly macroscopic state of the surrounding entities. This learning phase may involve both local interactions and statistical information extracted from the distribution functions. Subsequently, a decision phase takes place, in which the entity updates its microscopic state according to probabilistic transition rules.

A detailed analytical investigation of the interaction mechanisms for a-particles within the KTAP framework, including the formal derivation of transition probabilities, interaction rates, and their structural properties, is developed in~\cite{[BBL26]}. The multiscale formulation proposed here builds upon that analytical structure, extending it to incorporate cross-scale coupling and regulatory feedback mechanisms.

In spatially homogeneous settings, the sensitivity domain may reduce to a dependence on the internal activity variable only. In spatially distributed systems, however, the geometry of $\Omega$ and its dependence on the velocity direction introduce anisotropic and nonlocal effects, which are essential for modeling collective behaviors such as alignment, aggregation, or competitive exclusion.

The combination of intra-scale, cross-scale, conservative, and non-conservative interactions defines the general multiscale kinetic architecture. The corresponding balance equations are obtained by equating the transport of entities in the extended phase space to the net flux generated by all interaction mechanisms, thus leading to a coupled system of integro-differential equations for the distribution functions $f_i$ and $\varphi_j$.

\vskip.2cm \noindent \textit{\bf Interaction rates.}

The quantitative description of interactions requires the introduction of suitable \textit{interaction rates}, which measure the frequency with which encounters occur between entities belonging to the same or to different scales.

At the microscopic scale, interactions between a-particles belonging to the $i$-th and $h$-th functional subsystems are characterized by the interaction rate
\begin{equation}\label{alpha}
	\alpha_{ih}[f_i,f_h]
	(\bx,\bx^\ast,\bv_\ast,\bv^\ast,u_\ast,u^\ast),
\end{equation}
which measures the frequency of encounters between a candidate a-particle of the $i$-th subsystem, in the state $(\bx,\bv_\ast,u_\ast)$, and a field a-particle of the $h$-th subsystem, in the state $(\bx^\ast,\bv^\ast,u^\ast)$. The dependence on $f_i$ and $f_h$ accounts for the influence of the local statistical configuration of the system.

At the sub-microscopic scale, interactions between sa-particles belonging to the $j$-th and $k$-th sub-functional subsystems are characterized by
\begin{equation}\label{beta}
	\beta_{jk}[\varphi_j,\varphi_k]
	(\bx,\bx^\ast,\bv_\ast,\bv^\ast,w_\ast,w^\ast),
\end{equation}
which measures the frequency of encounters between a candidate sa-particle of the $j$-th subsystem and a field sa-particle of the $k$-th subsystem.

Cross-scale interactions are described by two distinct families of rates. Interactions in which an a-particle of the $i$-th FS interacts with an sa-particle of the $k$-th SFS are characterized by
\begin{equation}\label{gamma_ik}
	\gamma_{ik}[f_i,\varphi_k]
	(\bx,\bx^\ast,\bv_\ast,\bv^\ast,u_\ast,w^\ast),
\end{equation}
where we recall that the field entity is an sa-particle described by $(\bx^\ast,\bv^\ast,w^\ast)$. Conversely, interactions in which an sa-particle of the $j$-th SFS interacts with an a-particle of the $h$-th FS are characterized by
\begin{equation}\label{gamma_jh}
	\gamma_{jh}[\varphi_j,f_h]
	(\bx,\bx^\ast,\bv_\ast,\bv^\ast,w_\ast,u^\ast),
\end{equation}
where now the field entity is an a-particle described by $(\bx^\ast,\bv^\ast,u^\ast)$.

The explicit dependence of these rates on microscopic states, spatial variables, and distribution functions reflects the fact that interactions in living systems are generally nonlocal, anisotropic, and sensitive to the surrounding statistical configuration. Their precise analytical structure depends on the modeling assumptions adopted in specific applications.

\vskip.2cm \noindent \textit{\bf Transition probabilities.}

The outcome of interactions is described by suitable \textit{transition probability densities}, which model the probability that a candidate entity acquires a new microscopic state after interacting with a field entity. These probabilities encode the decision-making phase that follows the learning process occurring within the sensitivity domain.

\vskip.2cm \noindent \textbf{Microscopic FS--FS interactions}:
Consider a test a-particle of the $i$-th functional subsystem located at $(\bx,\bv,u)$. A candidate a-particle of the same subsystem is characterized by $(\bx,\bv_\ast,u_\ast)$, while a field a-particle belonging to the $h$-th subsystem is characterized by $(\bx^\ast,\bv^\ast,u^\ast)$.

The transition probability density is denoted by
\begin{equation}\label{A_ih}
	A_{ih}[f_i,f_h]
	\big(
	\bv_\ast \to \bv,\;
	u_\ast \to u
	\mid
	\bx, \bx^\ast,
	\bv_\ast, \bv^\ast,
	u_\ast, u^\ast
	\big),
\end{equation}
and represents the probability density that the candidate a-particle acquires the state $(\bv,u)$ after interaction with the field a-particle. The dependence on $f_i$ and $f_h$ reflects the influence of the surrounding statistical configuration on the decision process.

\vskip.2cm \noindent \textbf{Sub-microscopic SFS--SFS interactions}:
Consider a test sa-particle of the $j$-th sub-functional subsystem located at $(\bx,\bv,w)$. A candidate sa-particle is characterized by $(\bx,\bv_\ast,w_\ast)$, while a field sa-particle belonging to the $k$-th subsystem is characterized by $(\bx^\ast,\bv^\ast,w^\ast)$.

The corresponding transition probability density is
\begin{equation}\label{C_jk}
	C_{jk}[\varphi_j,\varphi_k]
	\big(
	\bv_\ast \to \bv,\;
	w_\ast \to w
	\mid
	\bx, \bx^\ast,
	\bv_\ast, \bv^\ast,
	w_\ast, w^\ast
	\big),
\end{equation}
which gives the probability density that the candidate sa-particle acquires the state $(\bv,w)$ after interaction. Also in this case, the dependence on $\varphi_j$ and $\varphi_k$ accounts for the collective statistical effects at the sub-microscopic scale.

\vskip.2cm \noindent \textit{\bf Transition probabilities for cross interactions.}

For conservative cross interactions between functional subsystems (FSs) and sub-functional subsystems (SFSs), we introduce the following transition probability densities.

\vskip.2cm \noindent \textbf{FS--SFS interactions}:
Let a candidate a-particle of the $i$-th functional subsystem be characterized by $(\bx,\bv_\ast,u_\ast)$, and let a field sa-particle of the $k$-th sub-functional subsystem be characterized by $(\bx^\ast,\bv^\ast,w^\ast)$. The transition probability density is defined by
\begin{equation}\label{B_ik}
	B_{ik}[f_i,\varphi_k]
	\big(
	\bv_\ast \to \bv,\;
	u_\ast \to u
	\mid
	\bx, \bx^\ast,
	\bv_\ast, \bv^\ast,
	u_\ast, w^\ast
	\big).
\end{equation}

\vskip.2cm \noindent \textbf{SFS--FS interactions}:
Let a candidate sa-particle of the $j$-th sub-functional subsystem be characterized by $(\bx,\bv_\ast,w_\ast)$, and let a field a-particle of the $h$-th functional subsystem be characterized by $(\bx^\ast,\bv^\ast,u^\ast)$. The transition probability density is defined by
\begin{equation}\label{D_jh}
	D_{jh}[\varphi_j,f_h]
	\big(
	\bv_\ast \to \bv,\;
	w_\ast \to w
	\mid
	\bx, \bx^\ast,
	\bv_\ast, \bv^\ast,
	w_\ast, u^\ast
	\big).
\end{equation}

\vskip.2cm \noindent \textbf{Remark.}
Both test a-particles and sa-particles lose their microscopic state as a consequence of conservative interactions.

\vskip.2cm \noindent \textit{\bf General differential system for conservative interactions.}

Limiting the analysis to conservative interactions, the balance equations are obtained by equating the transport term to the net interaction flux in the elementary volume of the space of microscopic states.

For the distribution functions $f_i$ and $\varphi_j$ one obtains the coupled system
\begin{equation}\label{coupled_system}
	\begin{cases}
		\displaystyle
		\left(
		\frac{\partial}{\partial t} + \bv \cdot \nabla_{\bx} \right) f_i(t,\bx,\bv,u) =
		\displaystyle \sum_{h=1}^{n} A_{ih}[f_i,f_h](t,\bx,\bv,u)\\[4mm]
		\displaystyle \hskip4cm + \sum_{k=1}^{m} B_{ik}[f_i,\varphi_k](t,\bx,\bv,u),\\[5mm]
		\displaystyle
		\left(
		\frac{\partial}{\partial t} + \bv \cdot \nabla_{\bx}
		\right) \varphi_j(t,\bx,\bv,w) =
		\displaystyle \sum_{k=1}^{m} C_{jk}[\varphi_j,\varphi_k](t,\bx,\bv,w)\\[5mm]
		\displaystyle \hskip4cm  + \sum_{h=1}^{n} D_{jh}[\varphi_j,f_h](t,\bx,\bv,w).
	\end{cases}
\end{equation}

Let $\Omega$ and $\Omega_m$ denote the sensitivity domains of a-particles and sa-particles, respectively. Let $D_{\bv}$ denote the velocity domain, and let $D_u$ and $D_w$ denote the activity domains of a-particles and sa-particles.

Define the integration domains
\begin{equation}\label{integration_domains_gamma}
	\Gamma = \Omega \times D_{\bv}^2 \times D_u^2,
	\qquad
	\Gamma_m = \Omega_m \times D_{\bv}^2 \times D_w^2,
\end{equation}
\begin{equation}\label{integration_domains_phi}
	\Phi = \Omega \times D_{\bv}^2 \times D_u \times D_w,
	\qquad
	\Phi_w = \Omega_m \times D_{\bv}^2 \times D_w \times D_u.
\end{equation}

\vskip.2cm \noindent \textit{\bf Explicit structure of microscopic operators.}

The conservative interaction terms acting on $f_i$ are given by gain--loss expressions.

\vskip.2cm \noindent \textbf{Microscopic interactions (FS--FS):}
\begin{equation}\label{micro_interactions_FS}
	\begin{aligned}
		A_i
		&=
		\sum_{h=1}^{n} \int_{\Gamma} \alpha_{ih}[f_i,f_h] (\bx,\bx^\ast,\bv_\ast,\bv^\ast,u_\ast,u^\ast) \\[2mm]
		& \qquad \times A_{ih}[f_i,f_h] (\bv_\ast \to \bv,\; u_\ast \to u
		\mid \bx,\bx^\ast,\bv_\ast,\bv^\ast,u_\ast,u^\ast) \\[2mm]
		&\qquad \times
		f_i(t,\bx,\bv_\ast,u_\ast)
		f_h(t,\bx^\ast,\bv^\ast,u^\ast)
		\, d\bx^\ast d\bv_\ast d\bv^\ast du_\ast du^\ast \\[3mm]
		& - f_i(t,\bx,\bv,u)
		\sum_{h=1}^{n}
		\int_{\Omega \times D_{\bv}\times D_u}
		\alpha_{ih}[f_i,f_h]
		(\bx,\bx^\ast,\bv,\bv^\ast,u,u^\ast) \\[2mm]
		&\qquad \times
		f_h(t,\bx^\ast,\bv^\ast,u^\ast)
		\, d\bx^\ast d\bv^\ast du^\ast .
	\end{aligned}
\end{equation}

\vskip.2cm \noindent \textbf{Cross-scale interactions (FS--SFS):}
\begin{equation}\label{cross_interactions_FS}
	\begin{aligned}
		B_i
		&= \sum_{k=1}^{m} \int_{\Phi}
		\gamma_{ik}[f_i,\varphi_k]
		(\bx,\bx^\ast,\bv_\ast,\bv^\ast,u_\ast,w^\ast) \\[2mm]
		&\qquad \times
		B_{ik}[f_i,\varphi_k]
		(\bv_\ast \to \bv,\;
		u_\ast \to u
		\mid
		\bx,\bx^\ast,\bv_\ast,\bv^\ast,u_\ast,w^\ast) \\[2mm]
		&\qquad \times
		f_i(t,\bx,\bv_\ast,u_\ast)
		\varphi_k(t,\bx^\ast,\bv^\ast,w^\ast)
		\, d\bx^\ast d\bv_\ast d\bv^\ast du_\ast dw^\ast \\[3mm]
		&
		- f_i(t,\bx,\bv,u)
		\sum_{k=1}^{m}
		\int_{\Omega \times D_{\bv}\times D_w}
		\gamma_{ik}[f_i,\varphi_k]
		(\bx,\bx^\ast,\bv,\bv^\ast,u,w^\ast) \\[2mm]
		&\qquad \times
		\varphi_k(t,\bx^\ast,\bv^\ast,w^\ast)
		\, d\bx^\ast d\bv^\ast dw^\ast .
	\end{aligned}
\end{equation}

\vskip.2cm \noindent \textit{\bf Explicit structure of sub-microscopic operators}

\vskip.2cm \noindent \textbf{Sub-microscopic SFS--SF interactions}:
\begin{equation}\label{sub-micro_interactions}
	\begin{aligned}
		C_j
		&=
		\sum_{k=1}^{m}
		\int_{\Gamma_m}
		\beta_{jk}[\varphi_j,\varphi_k]
		(\bx,\bx^\ast,\bv_\ast,\bv^\ast,w_\ast,w^\ast) \\[2mm]
		&\qquad \times
		C_{jk}[\varphi_j,\varphi_k]
		(\bv_\ast \to \bv,\;
		w_\ast \to w
		\mid
		\bx,\bx^\ast,\bv_\ast,\bv^\ast,w_\ast,w^\ast) \\[2mm]
		&\qquad \times
		\varphi_j(t,\bx,\bv_\ast,w_\ast)
		\varphi_k(t,\bx^\ast,\bv^\ast,w^\ast)
		\, d\bx^\ast d\bv_\ast d\bv^\ast dw_\ast dw^\ast \\[3mm]
		&
		- \varphi_j(t,\bx,\bv,w)
		\sum_{k=1}^{m}
		\int_{\Omega_m \times D_{\bv}\times D_w}
		\beta_{jk}[\varphi_j,\varphi_k]
		(\bx,\bx^\ast,\bv,\bv^\ast,w,w^\ast) \\[2mm]
		&\qquad \times
		\varphi_k(t,\bx^\ast,\bv^\ast,w^\ast)
		\, d\bx^\ast d\bv^\ast dw^\ast .
	\end{aligned}
\end{equation}

\vskip.2cm \noindent \textbf{Cross-scale SFS--FS interactions}:
\begin{equation}\label{cross_interactions_SFS}
	\begin{aligned}
		D_j
		&=
		\sum_{h=1}^{n}
		\int_{\Phi_w}
		\gamma_{jh}[\varphi_j,f_h]
		(\bx,\bx^\ast,\bv_\ast,\bv^\ast,w_\ast,u^\ast) \\[2mm]
		&\qquad \times
		D_{jh}[\varphi_j,f_h]
		(\bv_\ast \to \bv,\;
		w_\ast \to w
		\mid
		\bx,\bx^\ast,\bv_\ast,\bv^\ast,w_\ast,u^\ast) \\[2mm]
		&\qquad \times
		\varphi_j(t,\bx,\bv_\ast,w_\ast)
		f_h(t,\bx^\ast,\bv^\ast,u^\ast)
		\, d\bx^\ast d\bv_\ast d\bv^\ast dw_\ast du^\ast \\[3mm]
		&
		- \varphi_j(t,\bx,\bv,w)
		\sum_{h=1}^{n}
		\int_{\Omega_m \times D_{\bv}\times D_u}
		\gamma_{jh}[\varphi_j,f_h]
		(\bx,\bx^\ast,\bv,\bv^\ast,w,u^\ast) \\[2mm]
		&\qquad \times
		f_h(t,\bx^\ast,\bv^\ast,u^\ast)
		\, d\bx^\ast d\bv^\ast du^\ast .
	\end{aligned}
\end{equation}
\vskip.2cm \noindent \textit{\bf Proliferative and destructive interactions.}

Number-conservative interactions modify the microscopic state of the entities without altering their total number. However, a fundamental structural feature of living systems is the presence of proliferative and destructive events, which generate or remove entities as a consequence of interactions.

These mechanisms are intrinsically non-conservative and may act both at the microscopic (FS) and sub-microscopic (SFS) scales. They are typically triggered by cross-functional interactions and model birth, death, mutation, activation, or suppression phenomena.

\vskip.2cm \noindent \textbf{Microscopic proliferative and destructive operators}:
Let a test a-particle of the $i$-th FS interact with field a-particles of the $h$-th FS. We introduce the proliferative and destructive kernels
\begin{equation}\label{kernels}
	P_{ih}[f_i,f_h]
	(\bx,\bx^\ast,\bv,\bv^\ast,u,u^\ast),
	\qquad
	L_{ih}[f_i,f_h]
	(\bx,\bx^\ast,\bv,\bv^\ast,u,u^\ast),
\end{equation}
which respectively represent the rate of production and removal of test a-particles in the state $(\bx,\bv,u)$ induced by interactions with field particles in the state $(\bx^\ast,\bv^\ast,u^\ast)$.

The corresponding non-conservative operator acting on $f_i$ is written as
\begin{equation}\label{non-conservative_operator}
	\begin{aligned}
		E_i
		&=
		\sum_{h=1}^{n}
		\int_{\Omega \times D_{\bv}\times D_u}
		\alpha_{ih}[f_i,f_h]
		(\bx,\bx^\ast,\bv,\bv^\ast,u,u^\ast) \\
		&\qquad \times
		\Big[
		P_{ih}[f_i,f_h]
		-
		L_{ih}[f_i,f_h]
		\Big]
		(\bx,\bx^\ast,\bv,\bv^\ast,u,u^\ast) \\
		&\qquad \times
		f_i(t,\bx,\bv,u)
		\, f_h(t,\bx^\ast,\bv^\ast,u^\ast)
		\, d\bx^\ast d\bv^\ast du^\ast .
	\end{aligned}
\end{equation}

\vskip.2cm \noindent \textbf{Sub-microscopic proliferative and destructive operators}:
Analogously, proliferative and destructive events may affect sa-particles through interactions with a-particles. We introduce the kernels
\begin{equation}\label{m_kernels}
	P^{m}_{jh}[\varphi_j,f_h],
	\qquad
	L^{m}_{jh}[\varphi_j,f_h],
\end{equation}
modeling gain and loss of the $j$-th SFS induced by interaction with a-particles of the $h$-th FS.

The corresponding operator acting on $\varphi_j$ reads
\begin{equation}\label{F_j}
	\begin{aligned}
		F_j
		&=
		\sum_{h=1}^{n}
		\int_{\Omega_m \times D_{\bv}\times D_u}
		\gamma_{jh}[\varphi_j,f_h]
		(\bx,\bx^\ast,\bv,\bv^\ast,w,u^\ast) \\
		&\qquad \times
		\Big[
		P^{m}_{jh}[\varphi_j,f_h]
		-
		L^{m}_{jh}[\varphi_j,f_h]
		\Big]
		(\bx,\bx^\ast,\bv,\bv^\ast,w,u^\ast) \\
		&\qquad \times
		\varphi_j(t,\bx,\bv,w)
		\, f_h(t,\bx^\ast,\bv^\ast,u^\ast)
		\, d\bx^\ast d\bv^\ast du^\ast .
	\end{aligned}
\end{equation}

\vskip.2cm \noindent \textit{\bf General multiscale kinetic system.}

By combining conservative and non-conservative contributions, the complete multiscale kinetic system reads
\begin{equation}\label{complete_system}
	\begin{cases}
		\displaystyle
		\left(
		\frac{\partial}{\partial t}
		+
		\bv \cdot \nabla_{\bx}
		\right)
		f_i(t,\bx,\bv,u)
		=
		\bigg(A_i + B_i + E_i\bigg)[f, \varphi](t,\bx,\bv,u),
		\\[5mm]
		\displaystyle
		\left(
		\frac{\partial}{\partial t}
		+
		\bv \cdot \nabla_{\bx}
		\right)
		\varphi_j(t,\bx,\bv,w)
		=
		\bigg(C_j + D_j + F_j\bigg)[f, \varphi](t,\bx,\bv,w).
	\end{cases}
\end{equation}

A particularization of this differential system refers to the case of spatial homogeneity, which ca be written as follows:
\begin{equation}\label{spacehomogen}
	\begin{cases}
		\displaystyle \frac{\partial}{\partial t} f_i(t, u) = \bigg(A_i + B_i + E_i\bigg)[f, \vf](t, u),
		\\[5mm]
		\displaystyle \frac{\partial}{\partial t}\varphi_j(t,w) = \bigg(C_j + D_j + F_j\bigg)[f, \vf](t, w).
	\end{cases}
\end{equation}
Some applications related to this structure  will be discussed in Section 4.

\vskip.2cm \noindent \textbf{Remark 1.} The operators $E_i$ and $F_j$ break number conservation within each functional subsystem. In general, the total mass of the system is not preserved. Specific conservation laws may arise only under additional structural assumptions on the kernels $P_{ih}, L_{ih}, P^{m}_{jh}, L^{m}_{jh}$.

\vskip.2cm \noindent \textbf{Remark 2.} The mathematical architecture obtained in this section therefore provides a unified multiscale kinetic setting in which heterogeneity, adaptive decision processes, nonlinear interactions, and regulatory feedback mechanisms are consistently embedded. The next section will show how this general structure can be specialized into specific differential subsystems tailored to particular classes of living systems.

\section{From the theory to case studies and perspectives}\label{sec:4}

This section presents the basic guidelines for deriving multiscale models of the collective dynamics of real-world systems. According to the motivations given in the previous section, we first consider applications that require a detailed study of spatial dynamics. However, we also discuss how the case of spatial homogeneity can be considered in some specific applications and critically examine it. Therefore, we do not overlook spatial homogeneity; instead, we critically consider this aspect.

These topics are presented in the next two subsections, while a third subsection proposes a brief forward-looking perspective on further applications to be considered using the mathematical tools discussed in this paper. In this subsection, we specifically consider multi-dynamics systems characterized by interactions across different domains, for instance economics and politics, psychology and mechanics, and others presenting such features.

The presentation focuses solely on guidelines for deriving models. In this way, we address interested readers, encouraging them to develop further research programs accounting for the topics raised in this paper. Indeed, the aim of this section is to open a window onto prospective modeling applications based on the mathematical tools derived in the previous section.

\subsection{Multiscale tools for self-propelled particles}\label{subsec:4.1}

We apply the mathematical tools derived in the previous section to the multiscale modeling of the collective dynamics of self-propelled particles. More specifically, we refer to human crowds, swarms, programmed robots, and all systems that exhibit similar characteristics.

The literature on these systems contains interesting research focused on physical and mathematical modeling at the microscopic and mesoscopic levels when kinetic theory methods are applied (see~\cite{[BGO19],[BD06]}). Most of the literature focuses on a deterministic framework, in which the state of interacting entities is defined by mechanical and geometric variables. Our focus is on active particle methods, in which the microscopic level includes additional activity variables and sub-active particles play a role in collective dynamics.

For literature on self-propelled particles, we refer readers to the reviews on crowd dynamics methods in~\cite{[BLQRS23]} and on the dynamics of behavioral swarms in~\cite{[FLO26]}, which are formalized as systems of ordinary differential equations. Swarm theory refers to the deterministic model proposed by Felipe Cucker and Steven Smale (see~\cite{[CS07]}), who initiated research in this area, leading to a substantial body of literature. Recent applications have been proposed to model crowd dynamics using microscale models, as discussed by Liao and others (see~\cite{[BFLT26]}). In the absence of sa-particles, these papers have proposed models whose dynamics can be briefly described as follows:

\begin{enumerate}
	
	\vskip.1cm \item We consider a system of a-particles whose activity is the emotional state generated by the perception of danger.
	
	\vskip.1cm \item The micro- and macro-levels of the crowd in the sensitivity domain are learned by each a-particle.
	
	\vskip.1cm \item This learning is followed by a decision-making process that involves selecting the velocity direction and adapting the speed to the local density along the trajectory.
	
	\vskip.1cm \item The selection of velocity directions follows the above learning dynamics. First, a velocity direction is chosen by balancing the search for an exit, the avoidance of overcrowded areas, and the avoidance of close proximity to walls. This choice is weighted by the activity variable, the local mean velocity, and the density.
	
	\vskip.1cm \item Rational and irrational behaviors correspond to the middle level of activity, while only irrational behaviors correspond to high levels of activity. Irrational behavior involves imitating the mainstream, thus disregarding the search for less overcrowded areas.
	
\end{enumerate}

The guidelines for creating mathematical models at the mesoscopic scale (kinetic theory) are proposed in Section 4 of~\cite{[BGQR22]}, as summarized in the flowchart in Figure 5. Since this paper is open access, we do not repeat concepts that are already available elsewhere. However, we emphasize that avoiding overcrowded areas contributes to pedestrian safety.

The multiscale approach to accounting for vocal signals is useful for estimating the contribution of signals that encourage the avoidance of overcrowded areas to safety. Figure 1 shows the rationale for the multiscale approach.

\begin{figure}[t!]
	\begin{center}\scalebox{0.85}{
			\begin{tikzpicture}[node distance =3.0cm, auto]
				\node [block_long] (Physics) {\textbf{A}: Study \\of  behavioral \\ human crowds \\ with subsystems\\ of vocal signals};
				\node [block_long, below of=Physics,yshift=-1cm] (KTAP)  {\textbf{D}:  Derivation of \\ differential structures\\ modeling multiscale\\ collective dynamics};
				\node [block_long, left  of=KTAP,xshift=-2.0cm] (Micro) {\textbf{B}: Modeling SFSs: \\ dynamics and interactions\\ of vocal signals\\ within each SFSs and with FSs};
				\node [block_long, right of=KTAP,xshift= 2.0cm] (Macro)  {\textbf{C}: Modeling FSs \\ dynamics and interactions\\  mechanical and\\ activity variables};
				\node [block_long, below of=KTAP, yshift=-0.5cm,  fill=yellow!] (Learning) {\textbf{E}: Modelling \\ learning dynamics \\ within SFSs  and FSs\\ and across them};
				\node [block_long, below of=Learning,yshift=-0.5cm, fill=yellow!] (Decision) {\textbf{F}: Decision-making \\ $\to$  differential systems\\ describing the \\ collective dynamics};
				\draw [line width=.7mm, ->] (Physics.south) -- (Micro.north);
				\draw [line width=.7mm, ->] (Physics.south) -- (Macro.north);
				\draw [line width=.7mm, ->] (Physics.south) -- (KTAP.north);
				\draw [line width=.7mm, blue, ->] (KTAP.south) -- (Learning.north);
				\draw [line width=.7mm, blue, ->] (Micro.south) -- (Learning.west);
				\draw [line width=.7mm, blue, ->] (Macro.south) -- (Learning.east);
				\draw [line width=.7mm, blue, ->] (Learning.south) -- (Decision.north);
				\draw [line width=.7mm, ->] (KTAP.west) -- (Micro.east);
				\draw [line width=.7mm, ->] (KTAP.east) -- (Macro.west);
				\draw [line width=.7mm, ->] (Micro.east) -- (KTAP.west);
				\draw [line width=.7mm, ->] (Macro.west) -- (KTAP.east);
		\end{tikzpicture}}
	\end{center}
	\begin{center}
		\textit{Figure 1. On a rationale toward  multiscale modeling of behavioural crowds.}\label{fig2}
	\end{center}
\end{figure}

\subsection{On the dynamics in spatial homogeneity}\label{subsec:4.2}

This paper mainly focuses on collective dynamics, including time, space, velocity, and activity. The study specifically addresses living systems that exist far from equilibrium, as emphasized by Lee Hartwell, winner of the Nobel Prize in Physiology or Medicine, see~\cite{[HART99]}. According to this perspective, spatial homogeneity should be considered a special case in real-world living dynamics.

The modeling can be referred to the differential structure, see equation (\ref{spacehomogen}), where the microscopic dynamics refers to the activity variables of a-particles and s-particles, but does not take spatial dynamics into account. Some pioneering papers related to different aspects of immune competition have recently appeared (see~\cite{[BK26]} and~\cite{[BKL26]}). These papers describe the interaction between immune cells and disease carriers, which is mediated by sub-particles at the level of SFSs. These sub-particles then mediate the interaction-competition dynamics at the level of FSs. An additional study in~\cite{[BK24]} connects the contagion dynamics of individuals carrying infections to the in-host dynamics developed at the cellular level in epidemics.

Therefore, a systematic study of the multiscale kinetic theory of immune competition should be conducted within a general framework that can be related to different types of competition, including the impact that this type of study can have on medical care and society. The key challenge lies in developing models of biological systems that include spatial dynamics. This topic has been explored using active particles at the cellular level, see~\cite{[CKS21]} and~\cite{[CDS23]}. These studies could form part of a broader multiscale theory.

\subsection{Modeling and analytical perspectives}\label{subsec:4.3}

The previous subsections have critically analyzed some specific applications and indicated potential modeling perspectives. This concluding subsection briefly examines further perspectives. With a focus on applications, we emphasize that the aim of this paper is to cover various topics in the life sciences within an interdisciplinary framework.

For example, we apply the principles of behavioral economics, as outlined in~\cite{[BE24]}, to demonstrate that the kinetic theory of active particles can serve as a conceptual mathematical framework for Herbert A. Simon's artificial world, see~\cite{[Simon1962],[Simon2019]}. Further development of this research line, according to a multiscale approach that includes dynamics at lower scales corresponding to different types of individual communication, is worthwhile.

The derivation of macroscopic models from the underlying microscopic description is a research topic motivated by Hilbert's sixth problem, see~\cite{[Hilbert]}. There is an extensive body of literature on this problem for systems of classical particles, specifically focusing on the collective dynamics described by the Boltzmann equation. The interested reader is referred to the studies of Laure Saint-Raymond, see~\cite{[Laure]}, for the solution of Hilbert's problem developed as a perturbation of the Maxwellian equilibrium for this celebrated model of mathematical physics.

However, these studies cannot be simply extended to the kinetic theory of active particles. In fact, the differential system that describes the collective dynamics does not have a Maxwellian-type equilibrium. Results on this topic have instead been obtained by scaling the dynamics and considering small parameters as perturbations of the spatially homogeneous dynamics, see~\cite{[BC23]}.

The open problem is the micro–macro derivation of models including sub-microscopic dynamics. The challenge consists in exploring how the dynamics at the sub-microscopic scale would modify the models at the macroscopic scale.

\end{document}